\documentclass[12pt]{amsart}
\usepackage{amscd, amssymb}

\theoremstyle{plain}

\newtheorem{lemma}{Lemma}
\newtheorem*{mainlemma}{Main Lemma}
\newtheorem{theorem}{Theorem}
\newtheorem*{main}{Main Theorem}

\newtheorem*{claim*}{Claim}

\newtheorem*{question*}{Question}

\theoremstyle{definition}

\newtheorem*{definition*}{Definition}

\theoremstyle{remark}

\newtheorem*{remark*}{Remark}

\newtheorem*{acknowledgement}{Acknowledgement}

\newcommand{\cc}{\mathbb{C}}

\newcommand{\mbar}{{\overline{M}}}

\newcommand{\p}{\partial}

\newcommand{\ev}{{\operatorname{ev}}}

\newcommand{\vir}{\operatorname{vir}}
\newcommand{\la}{\langle}
\newcommand{\ra}{\rangle}

\newcommand{\h}{\hbar}

\newcommand{\bt}{t}

\newcommand{\HH}{\mathcal{H}}

\newcommand{\on}{\operatorname}

\renewcommand{\L}{\mathcal{L}}
\newcommand{\bq}{q}
\newcommand{\bp}{p}
\newcommand{\A}{\mathcal A}

\newcommand{\f}{f}
\newcommand{\g}{g}
\newcommand{\1}{\mathbf{1}}

\newcommand{\rspin}{\text{$r$-spin}}

\title[Invariance of tautological equations]
{Witten's conjecture, Virasoro conjecture,
and invariance of tautological equations}

\author{Y.-P.~Lee}

\address{Department of Mathematics \\
        University of Utah \\
        Salt Lake City, Utah 84112-0090\\
        U.S.A.}
\email{yplee@math.utah.edu}

\address{NCTS, Mathematics Division\\
    4F Third General Building, National Tsing Hua University\\
    No. 101, Sec. 2, Kuang Fu Road\\
    Hsinchu, Taiwan 30043}
\email{yplee@math.cts.nthu.edu.tw}

\thanks{Research partially supported by NSF}

\begin{document}

\begin{abstract}
The main goal of this paper is to prove the following two conjectures
\emph{for genus up to two}:
\begin{enumerate}
\item Witten's conjecture on the relations between higher spin curves and
Gelfand--Dickey hierarchy.
\item Virasoro conjecture for target manifolds with conformal semisimple
Frobenius manifolds.
\end{enumerate}
The main technique used in the proof is the invariance of tautological
equations under loop group action.
\end{abstract}

\maketitle

\setcounter{section}{-1}

\section{Introduction}

\subsection{Two dimensional quantum gravity}

The famous conjecture by E.~Witten \cite{EW1} in 1990 predicted a striking
relation between two seemingly unrelated objects: A generating functions of
intersection numbers on moduli spaces of stable curves and a $\tau$-function of
the KdV hierarchy. The physical basis of this conjecture comes from the
identification of two approaches to the two dimensional quantum gravity.
Roughly, the correlators of the two dimensional quantum gravity are Feynman
path integrals over the ``space of metrics'' on two dimensional topological
real surfaces. One approach of evaluating this path integral involves a
topological field theory technique which is expected to reduce to the
integration over the moduli space of curves. The other approach considers an
approximation of the space of the metrics by piecewise flat metrics and then
take a suitable continuous limit.

In the first approach, the free energy becomes the following geometrically
defined function
\[
 \tau^{pt}(t_0,t_1, \ldots)
 = e^{\sum_{g=0}^{\infty} \hbar^{g-1} F^{pt}_g(t_0,t_1, \ldots)} ,
\]
where $F^{pt}_g(t)$ is the generating function of (tautological) intersection
numbers on the moduli space of stable curves of genus $g$
\[
  F^{pt}_g (t_0,t_1, \ldots) := \sum_n \frac{1}{n!} \int_{\mbar_{g,n}}
  \prod_{i=1}^n (\sum_k t_k \psi^k_i) .
\]
($\hbar$ is usually set to be 1 in the literature.) Moreover, from elementary
geometry of moduli spaces, one easily deduces that $\tau^{pt}$ satisfies an
additional equation, called the \emph{string equation} (or puncture equation).
It is a basic fact in the theory of KdV (or in general KP) hierarchies that the
string equation uniquely determines one $\tau$-function for the KdV hierarchy
from all $\tau$-functions parameterized by Sato's grassmannian.

In the second approach, the generating function in the double scaling limit
yields the $\tau$-function $\tau^{qg}$ of the KdV hierarchy whose initial value
is $d^2/dx^2 + 2x$. From here Witten asserts that $\tau^{pt}$ must be equal to
$\tau^{qg}$ since there should be only one quantum gravity.

\subsection{Witten's conjecture on spin curves and Gelfand--Dickey hierarchies}

In 1991 Witten formulated a remarkable generalization of the above conjecture.
He argued that an analogous generating function $\tau^{\rspin}$ of the
intersection numbers on moduli spaces of $r$-spin curves should be identified
as a $\tau$-function of Gelfand--Dickey ($r$-KdV) hierarchies \cite{EW2}. When
$r=2$, this conjecture eventually reduces to the previous one as $2$-KdV is the
ordinary KdV.

The special case $\tau^{pt}=\tau^{qg}$ was soon proved by M.~Kontsevich
\cite{MK}. More recently a new proof was given by Okounkov--Pandharipande
\cite{OP}. However, the generalized conjecture remains open up to this day.

Throughout the 12 years, there has been substantial progress in the
foundational issues involved in the 1991 conjecture. In particular,
Jarvis--Kimura--Vaintrob \cite{JKV} established the genus zero case of the
conjecture; T.~Mochizuki and A.~Polishchuk independently established the
following property for $\tau^{\rspin}$:

\begin{theorem} \cite{TM, AP}
All tautological equations hold for $F_g^{\rspin}$.
\footnote{Tautological relations in this article means both the relations of
the tautological classes on moduli spaces of curves and the induced relations
in the ``cohomological field theories''.}
\end{theorem}
In fact, $F_g^{\rspin}$ satisfies all ``expected functorial properties'',
similar to the axioms formulated by Kontsevich--Manin in the Gromov--Witten
theory.

However, Riemann's trichotomy of Riemann surfaces has taught us that things are
very different in genus one and in genus $\ge 2$. Our Main Theorem therefore
provides a solid confirmation for Witten's 1991 conjecture, covering one
example ($g=1$ and $g=2$) for the other two cases of the trichotomy. \emph{In
fact, this work starts as a project trying to understand this conjecture in
higher genus.}

For more background information about Witten's conjecture, the readers are
referred to Witten's original article \cite{EW2} and the paper \cite{JKV} by
Jarvis--Kimura--Vaintrob.
In the remaining of this article, ``Witten's conjecture'' means the 1991
conjecture if not otherwise specified.

\subsection{Virasoro conjecture}

In 1997 another generalization of Witten's 1990 conjecture was proposed by
T.~Eguchi, K.~Hori and C.~Xiong. Witten's 1990 conjecture has an equivalent
formulation \cite{DVV} \cite{FKN} \cite{MK}: $\tau^{pt}$ is annihilated by
infinitely many differential operators $\{ L^{pt}_n \}, n \ge -1$, satisfying
the Virasoro relations
\[
  [L_m, L_n]=(m-n) L_{m+n}
\]
such that $L_{-1} \tau^{pt} =0$ is the string equation alluded above. To
generalize Witten's 1990 conjecture to any projective smooth variety $X$,
consider the moduli spaces of curves as the moduli spaces of maps to a point.
It is clear that $\tau^{pt}$ should be replaced by
\begin{equation} \label{e:1}
 \tau_{GW}^X (\bt) := e^{\sum_{g=0}^{\infty} \hbar^{g-1} F^X_g(\bt)},
\end{equation}
where $F^X_g (\bt)$ is the generating function of genus $g$ Gromov--Witten
invariants with descendents for $X$.
Based on that Eguchi--Hori--Xiong \cite{EHX}, and S.~Katz, managed to to define
$\{ L^X_n \}$ for $n \ge -1$, satisfying the Virasoro relations. \footnote{By
the Virasoro relations, one only has to construct $L_2^X$, and the rest will
follow.} They conjectured that
\[
  L^X_n \tau^X (\bt) =0, \quad \text{for $n \ge -1$}.
\]
This conjecture is commonly referred to as the \emph{Virasoro conjecture}.

Eguchi--Hori--Xiong gave a partial proof for their conjecture in genus zero and
a proof of $L_0^X \tau^X=0$, among other things. Later X.~Liu and G.~Tian
\cite{LT} proved the genus zero case in general. Using a very different method,
Dubrovin--Zhang \cite{DZ2} established the genus one case of Virasoro
conjecture for \emph{conformal} semisimple Frobenius manifolds. \footnote{The
definition of Frobenius manifolds in this article does not require existence of
an Euler field, which is assumed in Dubrovin's definition. Dubrovin's
definition will be referred to as \emph{conformal Frobenius manifold} instead.}
The recent progress by Givental \cite{AG2} and by Okounkov--Pandharipande
\cite{OP} have confirmed the conjecture for toric Fano manifolds and curves
respectively at all genus. Some background information on Virasoro conjecture
can be found in \cite{EG3} and \cite{LP}.

\subsection{Main results}

\begin{main}
Witten's conjecture and Virasoro conjecture for manifolds with conformal
semisimple quantum cohomology hold up to genus two.
\end{main}

A few words about the main idea in the proof. Givental in a series of papers
\cite{AG1} \cite{AG2} \cite{AG3} \cite{AG4} introduces a definition of higher
genus potentials for any semisimple Frobenius manifold which is not necessarily
the quantum cohomology of a projective manifold. This definition is
``formulaic'' in the sense that the higher genus potentials are \emph{defined}
by a formula from the data of semisimple Frobenius manifolds (i.e.~genus zero
data). This enables him to prove that his theory satisfies Virasoro conjecture
and, in the case of $A_n$ singularities, Witten's conjecture. However,
Givental's theory is conjecturally equivalent to the geometric theory, whether
it is the Gromov--Witten theory or the theory of spin curves. Therefore, what
is needed here is a proof that Givental's theory is equal to the geometric
theory.

The bulk of this paper is devoted to this proof at genus two. Similar statement
in genus one is proved in the conformal case in \cite{DZ1} and in the general
case in \cite{GL}.

\begin{remark*}
There are other possible approaches to this problem. Our earlier approach in
\cite{YL1} reduces the checking of the Main Theorem to complicated, but
finite-time checkable, identities. Nevertheless, it lacks the underlying
simplicity of this approach.

After this result was announced, X.~Liu \cite{XL2} informed us that he was also
able to reduce the genus two Virasoro conjecture for semisimple Gromov--Witten
theory to some complicated identities and he was able to check these identities
by hand and by a Mathematica program.
\end{remark*}

\begin{acknowledgement}
This idea of this work first comes in the form of \cite{YL1} while working
jointly with R.~Pandharipande in the book project \cite{LP}, and this current
approach has its root in a recent joint work with A.~Givental \cite{GL}. It is
a great pleasure to thank both of them. Thanks are also due to T.~Jarvis,
T.~Kimura, X.~Liu, Y.~Ruan, A.~Vaintrob and especially E.~Getzler for useful
discussions and communications. Part of this work was done during a visit to
NCTS, whose hospitality is greatly appreciated.
\end{acknowledgement}

\section{Review of geometric Gromov--Witten theory}

Gromov--Witten theory studies the tautological intersection theory on
$\mbar_{g,n}(X,\beta)$, the moduli spaces of stable maps from curves $C$ of
genus $g$ with $n$ marked points to a smooth projective variety $X$. The
intersection numbers, or \emph{Gromov--Witten invariants}, are integrals of
tautological classes over the virtual fundamental classes of
$\mbar_{g,n}(X,\beta)$
\[
  \int_{[\mbar_{g,n}(X,\beta)]^{\vir}}
    \prod_{i=1}^n \ev_i^*(\gamma_i) \psi_i^{k_i}.
\]
Here $\gamma_i \in H^*(X)$
and $\psi_i$ are the \emph{cotangent classes} (gravitational
descendents).

For the sake of the later reference, let us fix some notations.
\begin{enumerate}
\item[(i)] $H := H^*(X, \cc)$, assumed of rank $N$.

\item[(ii)] Let $\{ \phi_{\mu} \}_{\mu=1}^N$ be an \emph{orthonormal} basis of
$H$ with respect to the Poincar\'e pairing.

\item[(iii)] Let $\mathcal{H}_t := \oplus_0^{\infty} H$ be the infinite
dimensional complex vector space with basis $\{ \phi_{\mu} \psi^k \}$.

\item[(iv)] Let $\{ t^{\mu}_k \}$, $\mu=1, \ldots, N$, $k=0, \ldots, \infty$,
be the dual coordinates of the basis $\{ \phi_{\mu} \psi^k \}$.

\end{enumerate}
We note that at each marked point, the insertion is $\mathcal{H}_t$-valued. Let
$t:= \sum_{k, \mu} t^{\mu}_k \phi_{\mu} \psi^k$ denote a general element in the
vector space $\mathcal{H}_t$.

\begin{enumerate}

\item[(v)] Define {$\displaystyle \la \p^{\mu_1}_{k_1} \ldots \p^{\mu_n}_{k_n}
\ra_{g,n,\beta} := \int_{[\mbar_{g,n}(X,\beta)]^{\vir}} \prod_{i=1}^n
\ev_i^*(\phi_{\mu_i}) \psi_i^{k_i}$} and define $\la t^n \ra_{g,n,\beta}=\la t
\ldots t \ra_{g,n,\beta}$ by multi-linearity.

\item[(vi)] Let {$\displaystyle F^X_g(t) := \sum_{n, \beta} \frac{1}{n!} \la
t^n \ra_{g,n,\beta}$} be the generating function of all genus $g$
Gromov--Witten invariants.

\end{enumerate}

The ``$\tau$-function of $X$'' is the formal expression {$\displaystyle
\tau_{GW}^X := e^{\sum_{g=0}^{\infty} \hbar^{g-1} F_g^X}$} defined in
\eqref{e:1}.

\subsection{Tautological equations} \label{s:1.1}
Let $E=0$ be a \emph{tautological equation}, i.e.~a equation of the
tautological classes in the moduli space of stable curves $\mbar_{g,n}$. Since
there is a morphism
\[
  \mbar_{g,n}(X,\beta) \to \mbar_{g,n}
\]
by forgetting the map, one can pull-back $E=0$ to $\mbar_{g,n}(X,\beta)$. Due
to the functorial properties of the virtual fundamental classes, the pull-backs
of the \emph{tautological equations hold for the Gromov--Witten theory of any
target space}. The term \emph{tautological equations} will also be used for the
corresponding equations in the Gromov--Witten theory and in the theory of spin
curves.

\section{Genus zero axiomatic theory}

Let $H$ be a complex vector space of dimension $N$ with a distinguished element
$\1$. Let $(\cdot,\cdot)$ be a $\cc$-bilinear metric on $H$, i.e.~a
nondegenerate symmetric $\cc$-bilinear form. Let $\HH$ denote the infinite
dimensional complex vector space $H((z^{-1}))$ consisting of Laurent formal
series in $1/z$ with vector coefficients. \footnote{Different completions of
this spaces are used in different places, but this will be not be discussed in
details in the present article as it is involving too much. See \cite{LP} for
details.} Introduce the symplectic form $\Omega$ on $\HH$:
\[ \Omega (\f,\g )  = \frac{1}{2\pi i} \oint\ ( \f (-z), \g(z) ) \ dz .\]
There is a natural polarization $\HH = \HH_{+}\oplus \HH_{-}$ by the Lagrangian
subspaces $\HH_{+} = H[z]$ and $\HH_{-} = z^{-1} H [[ z^{-1} ]]$  which
provides a symplectic identification of $(\HH, \Omega)$ with the cotangent
bundle $T^*\HH_{+}$.

Let $\{ \phi_{\mu} \}$ be an \emph{orthonormal} basis of $H$. An $H$-valued
Laurent formal series can be written in this basis as
\begin{multline*}
 \ldots + (\bp^1_1,\ldots,\bp^N_{1}) \frac{1}{(-z)^2}
 + (\bp^1_{0},\ldots, \bp^N_{0}) \frac{1}{(-z)} \\
 + (\bq^1_{0},\ldots, \bq^N_{0})
 + (\bq^1_{1}, \ldots, \bq^N_{1}) z + \ldots.
\end{multline*}
In fact, $\{ \bp_k^{\mu}, \bq_k^{\mu} \}$ for $k= 0, 1, 2, \ldots$ and
$\mu=1, \ldots, N$ are the Darboux coordinates compatible with this
polarization in the sense that
\[
 \Omega = \sum d \bp^{\mu}_k \wedge d \bq^{\mu}_k .
\]
To simplify the notations, $\bp_k$ will stand for the vector
$(\bp^1_{k},\ldots,\bp^N_{k})$ and $\bp^{\mu}$ for
$(\bp^{\mu}_0, \bp^{\mu}_1, \ldots )$.

The parallel between $\HH_+$ and $\mathcal{H}_t$ is evident, and is in fact
given by the affine coordinate transformation, the \emph{dilaton shift},
\[
  t^{\mu}_k = q^{\mu}_k + \delta^{\mu 1} \delta_{k 1}.
\]

\begin{definition*}
Let $G_0(t)$ be a (formal) function on $\HH_+$. The pair $(\HH, G_0)$ is called
a \emph{$g=0$ axiomatic theory} if $G_0$ satisfies three sets of genus zero
tautological equations: the \emph{Topological Recursion Relations} (TRR)
\eqref{e:g0trr}, the {\em String Equation} \eqref{e:se} and the {\em Dilaton
Equation} \eqref{e:de}.

\begin{align}
 \label{e:de}
 &\frac{\p F_0(\bt)}{\p \bt^1_1} (\bt)
 = \sum_{n=0}^{\infty} \sum_{\nu} \bt^{\nu}_n
 \frac{\p F_0(\bt)}{\p \bt^{\nu}_n} - 2F_0 (\bt), \\
  \label{e:se}
  &\frac{\p F_0 (\bt)}{\p \bt^1_0} =
  \frac{1}{2}(\bt_0,\bt_0)+ \sum_{n=0}^{\infty}
  \sum_{\nu} \bt_{n+1}^{\nu} \frac{\p F_0 (\bt)}{\p \bt^{\nu}_n},  \\
  \label{e:g0trr}
  &\frac{\p^3 F_0 (\bt)}
       {\p \bt^{\alpha}_{k+1} \p \bt^{\beta}_{l} \p \bt^{\gamma}_m}
  = \sum_{\mu} \frac{\p^2 F_0 (\bt)}{\p \bt^{\alpha}_{k} \p \bt^{\mu}_{0}}
    \frac{\p^3 F_0 (\bt)}
       {\p \bt^{\mu}_{0} \p \bt^{\beta}_{l} \p \bt^{\gamma}_m}
\end{align}
\end{definition*}

It can be shown that this is equivalent to the definition of abstract formal
Frobenius manifolds, not necessarily conformal. The coordinates on the
corresponding Frobenius manifold is given by the following map \cite{DW}
\begin{equation} \label{e:dg}
 s^{\mu} := \frac{\p}{\p \bt^{\mu}_0}\frac{\p}{\p \bt^1_0} G_0 (\bt).
\end{equation}

In the case of the geometric theory, one may elect to use either formulation.
The main advantage, seems to us, is the expansion of the formulation from
$\mathcal{H}_t$ to $\HH$ where a symplectic structure is available. In the
latter case, many properties can be reformulated in terms of the symplectic
structure $\Omega$ and hence independent of the choice of the polarization.
This suggests that the space of ``genus zero axiomatic Gromov--Witten
theories'', i.e. the space of $G_0$ satisfying the string equation, dilaton
equation, topological recursion relations (TRR), has a huge symmetry group.

\begin{definition*}
Let $\L^{(2)}GL(H)$ denote the {\em twisted loop group} which consists of
$\operatorname{End}(H)$-valued formal Laurent formal series $M(z)$ in
the indeterminate $z^{-1}$ satisfying $M^*(-z)M(z)=\1$.
Here $\ ^*$ denotes the adjoint with respect to $(\cdot ,\cdot )$.
\end{definition*}

The condition $M^*(-z)M(z)=\1$ means that $M(z)$ is a symplectic transformation
on $\HH$.

\begin{theorem} \label{t:2} \cite{AG4}
The twisted loop group acts on the space of axiomatic genus zero theories.
Furthermore, the action is transitive on the semisimple theories of a fixed
rank $N$.
\end{theorem}

\section{Quantization and higher genus potentials} \label{s:3}

\subsection{Preliminaries on quantization}

To quantize an infinitesimal symplectic transformation, or its corresponding
quadratic hamiltonians, we recall the standard Weyl quantization. A
polarization $\HH=T^* \HH_+$ on the symplectic vector space $\HH$ (the phase
space) defines a configuration space $\HH_+$. The quantum ``Fock space'' will
be a certain class of functions $f(\h, q)$ on $\HH_+$ (containing at least
polynomial functions), with additional formal variable $\hbar$ (``Planck's
constant''). The classical observables are certain functions of $\bp, \bq$. The
quantization process is to find for the classical mechanical system on $\HH$ a
``quantum mechanical'' system on the Fock space such that the classical
observables, like the hamiltonians $h(q,p)$ on $\HH$, are quantized to become
operators $\hat{h}(q,\dfrac{\p}{\p q})$ on the Fock space.

Let $A(z)$ be an $\on{End}(H)$-valued Laurent formal series in $z$ satisfying
\[
  (A(-z) f(-z), g(z)) + (f(-z), A(z) g(z)) =0,
\]
then $A(z)$ defines an infinitesimal symplectic transformation
\[
  \Omega(A f, g) + \Omega(f, A g)=0.
\]
An infinitesimal symplectic transformation $A$ of $\HH$ corresponds to a
quadratic polynomial $P(A)$ in $\bp, \bq$
\[
  P(A)(f) := \frac{1}{2} \Omega(Af, f) .
\]
For example, let $\dim H=1$ and $A(z) =1/z$. It is easy to see that $A(z)$ is
infinitesimally symplectic and
\begin{equation} \label{e:cse}
 P(z^{-1})= -\frac{\bq_0^2}{2} - \sum_{m=0}^{\infty} \bq_{m+1} \bp_m.
\end{equation}

In the above Darboux coordinates, the quantization $P \mapsto \hat{P}$ assigns
\begin{equation} \label{e:wq}
 \begin{split}
  &\hat{1}= 1, \  \hat{p}_k^i= \sqrt{\hbar} \frac{\p}{\p q^i_k},
   \hat{q}^i_k = q^i_k / {\sqrt{\hbar}}, \\
  &(p^i_k p^j_l) \hat{\ } = \hat{p}^i_k \hat{p}^j_l
    =\hbar \frac{\p}{\p q^i_k} \frac{\p}{\p q^j_l}, \\
   &(p^i_k q^j_l) \hat{\ } = q^j_l \frac{\p}{\p q^i_k},\\
  &(q^i_k q^j_l) \hat{\ } = \hat{q}^i_k \hat{q}^j_l /\hbar ,
 \end{split}
\end{equation}

Note that one often has to quantize the symplectic instead of the infinitesimal
symplectic transformations. Following the common practice in physics, we define
\begin{equation} \label{e:q}
  (e^{A(z)})\hat{}\ := e^{(A(z)) \hat{}\ } ,
\end{equation}
for $e^{A(z)}$ an element in the twisted loop group.

\subsection{$\tau$-function for semisimple Frobenius manifolds}

Let $H^{N pt}$ be the rank $N$ Frobenius manifold corresponding to $X$ being
$N$ points. In this case, the delta-functions at the $N$ points form an
orthonormal basis $\{ \phi_{\mu} \}$ and the idempotents of the quantum product
\[
  \phi_{\mu} * \phi_{\nu} = \delta_{\mu \nu} \phi_{\mu}.
\]
The genus zero potential is nothing but a product of genus zero potentials of
$N$ points
\[
  F^{N pt}_0 (\bt^1, \ldots, \bt^N) = F^{pt}_0 (\bt^1) \ldots F^{pt}_0(\bt^N).
\]
Note that $G_0^{H^{N pt}} = F_0^{N pt}$. By Theorem~\ref{t:2}, the genus zero
potential $G_0^H$ of any semisimple formal Frobenius manifold $H$ can be
obtained from $G_0^{H^{N pt}}$ by the action of an element $O_H$ in the twisted
loop group. By Birkhoff factorization, $O_H = S_H (z^{-1}) R_H (z)$, where
$S(z^{-1})$ (resp.~$R(z)$) is an matrix-valued functions in $z^{-1}$
(resp.~$z$). \footnote{In fact $R(z)$ is a series in $z$ and therefore not
really an element in the twisted loop group, but rather in its suitable
completion. See \cite{LP}.}

In order to define the higher genus potentials $G_g^H$, one first introduces
the ``$\tau$-function of $H$''
\begin{equation} \label{e:taug}
  \tau_G^H := \hat{S}_H \hat{R}_H \tau^{N pt},
\end{equation}
and define $G_g^H$ via the formula (cf.~\eqref{e:1})
\begin{equation} \label{e:2}
  \tau_G^H =: e^{\sum_{g=0}^{\infty} \hbar^{g-1} G_g^H}.
\end{equation}
Strictly speaking, the multiplication, $\hat{S}_H \hat{R}_H$, is not
well-defined. However, the function $\hat{S}_H (\hat{R}_H \tau^{N pt})$ is,
thanks to the $(3g-2)$-jet properties. We will not discuss this subtle point
here but refer the interested readers to \cite{LP}.

What makes the above model especially attractive are the facts that
\begin{enumerate}

\item[(a)] It works for any semisimple Frobenius manifolds.

\item[(b)] It enjoys properties often complementary to the geometric theory.
\end{enumerate}
Thanks to (a), one also has a definition for the Frobenius manifolds
$H_{A_{r-1}}$ of the miniversal deformation space of $A_{r-1}$ singularity. It
turns out that this Frobenius manifold is isomorphic to the Frobenius manifold
defined by the genus zero potential of $r$-spin curves. Furthermore, Givental
has proved

\begin{theorem} \cite{AG3} \label{t:ag3}
$\tau_G^{H_{A_{r-1}}}$ is a $\tau$-function of $r$-KdV hierarchy.
\end{theorem}

As in the case of the ordinary KdV, it is easy to show that both
$\tau_G^{H_{A_{r-1}}}$ and $\tau^{\rspin}$ satisfy the additional string
equation. Therefore, in order to prove Witten's conjecture, one only has to
show $G_g^{H_{A_{r-1}}}= F^{\rspin}_g$.

As for (b), note for example that the Virasoro constraints for $\tau_G^H$
follow almost from the definition. As discussed in Section~\ref{s:1.1},
$\tau_{GW}^X$ satisfies the tautological equations due to some functorial
properties built in the definition of the Gromov--Witten theory. However, the
Virasoro constraints for $\tau_{GW}^X$ and tautological equations for
$\tau_G^H$ are highly nontrivial challenges. An obvious, and indeed very good,
strategy to resolve all the above questions at once is to answer the following
question: \footnote{This is basically Givental's conjecture \cite{AG2},
although we have included $r$-spin curves on the geometric side, which is
strictly speaking not a geometric Gromov--Witten theory.}

\begin{question*}
Is $G_g = F_g$ That is, does Givental's construction coincide with the
geometric one when both are available?
\end{question*}


\subsection{Tautological equations and uniqueness theorems}

Our approach to the Question is to show that $G_g$ satisfies enough geometric
properties of $F_g$ so that they have to be equal by some uniqueness theorems.
More specifically, the geometric properties we will utilize are the
tautological equations. For simplicity of language, let us call the \emph{genus
zero tautological equations} the following genus zero equations: topological
recursion relations (TRR), string equation and dilaton equation; the
\emph{genus one tautological equations} the following two equations: genus one
Getzler's equation \cite{EG0} and genus one TRR; the \emph{genus two
tautological equations} the set of 3 equations by Mumford \eqref{e:m}, Getzler
\cite{EG1}, and Belorousski--Pandharipande (BP) \cite{BP}.

In genus one, Dubrovin and Zhang \cite{DZ1} made the following important
observation of the uniqueness property.

\begin{lemma} \cite{DZ1} \label{l:dz}
Let $G_0(\bt)$ be the genus zero potential of a semisimple Frobenius manifold
$H$. Suppose that both pairs $(G_0(\bt), F_1 (\bt))$ and $(G_0(\bt),G_1(\bt)$
satisfy genus one Getzler's equation and topological recursion relations. Then
$F_1 - G_1$ is a linear combination of canonical coordinates. Furthermore, if
$H$ is conformal and both pairs satisfy the conformal equation, then $F_1 -
G_1$ is a constant.
\end{lemma}

The proof of this fact goes as follows. First, genus one TRR guarantees that
the descendent invariants are uniquely determined by primary invariants.
Second, genus one Getzler's equation, when written in canonical coordinates
$u^i$, is equal to $\frac{\p^2 F_1}{\p u^i \p u^j} = B_{ij}$ where $B_{ij}$
involves only genus zero invariants. Moreover, the conformal structure
determined by a linear vector field (Euler field), uniquely determines the
linear term. We will refer to this fact casually as ``The genus one potential
for a conformal Frobenius manifold is uniquely determined by genus one
tautological equations.'' In the same spirit, the genus two uniqueness theorem
of X.~Liu is formulated:

\begin{theorem} \cite{XL} \label{t:xl}
The genus two descendent potentials for any conformal semisimple Frobenius
manifolds are uniquely determined by genus two tautological equations.
\end{theorem}

The proof of Liu's theorem rests on some very complicated calculation. We note
that it is not known at this point whether this uniqueness theorem, or any
weaker version, holds for non-conformal semisimple Frobenius manifolds.

\begin{remark*}
There is another type of uniqueness theorem: Dubrovin and Zhang \cite{DZ3} have
proved that Virasoro conjecture plus $(3g-2)$-jet property uniquely determines
$\tau$-function for any semisimple Frobenius manifold. The $(3g-2$)-jet
property is proved by Getzler \cite{EG2} in the geometric Gromov--Witten theory
and by Givental \cite{AG2} in the context of semisimple Frobenius manifolds. It
is also expected to hold for the $\tau^{\rspin}$. Therefore, a proof of the
Virasoro conjecture for $\tau^X$ should also answer the above Question
positively.
\end{remark*}

Note that
\begin{itemize}

\item $\tau_G = \hat{S} \hat{R} \tau^{N pt}$ and $\tau^{N pt} (t^1, \ldots t^N)
= \prod_{i=1}^N \tau^{pt} (t^i)$.

\item $\tau^{pt}(\bt)$ satisfies all tautological equations.
\end{itemize}
It follows that

\begin{mainlemma}
In order to show that a set of tautological equations holds for $G_g$,
it suffices to show that it is invariant under arbitrary $\on{End}(H)$-valued
series $\hat{S}(z^{-1})$ and $\hat{R}(z)$.
\end{mainlemma}

This lemma is our main technical tool to prove the Main Theorem. In fact, in
order to prove the invariance of the tautological equations, it is enough to
prove the \emph{infinitesimal} invariance of the tautological equations. Before
we proceed, let us study more carefully the quantized twisted loop groups.

\section{Quantization of twisted loop groups} The twisted loop group is
generated by ``lower triangular subgroup'' and the ``upper triangular
subgroup''. The lower triangular subgroup consists of $\on{End}(H)$-valued
formal series $S(z^{-1})= e^{s(z^{-1})}$ in $z^{-1}$ satisfying $S^*(-z) S(z) =
\mathbf{1}$ or equivalently
\[
 s^*(-z^{-1}) + s(z^{-1}) =0.
\]
The upper triangular subgroup consists of the regular part of the twisted loop
groups $R(z)= e^{r(z)}$ satisfying $R^*(-z) R(z) = \mathbf{1}$ or equivalently
\begin{equation} \label{e:r}
  r^*(-z) + r(z) =0.
\end{equation}
In fact, we will use $R(z)$ to denote an $\on{End}(H)$-valued series in $z$,
and call it an element in the ``upper triangular subgroup'' by abusing the
language.

\subsection{Quantization of lower triangular subgroups}

The quadratic hamiltonian of $s(z^{-1}) = \sum_{l=1}^{\infty} s_l z^{-l}$ is
\[
  \sum_{l=1}^{\infty} \sum_{n=0}^{\infty} \sum_{i,j}
  (s_{l})_{ij} q^j_{l+n} p^i_n +
  \sum \frac{1}{2} (-1)^n (s_{l})_{ij} q^i_n q^j_{l-n-1}.
\]
The fact that $s(z^{-1})$ is a series in $z^{-1}$ implies that the
quadratic hamiltonian $P(s)$ of $s$ is of
the form $q^2\text{-term} + qp\text{-term}$ where $q$ in $qp$-term
does not contain $q_0$. The quantization of the $P(s)$
\[
 \hat{s} = \sum (s_{l})_{ij} q^j_{l+n} \p_{q^i_n} +
 \frac{1}{2\hbar} \sum (-1)^{n} (s_{l})_{ij} q^i_n q^j_{l-n-1} .
\]
Here $i,j$ are the indices of the orthonormal basis. (The indices $\mu, \nu$
will be reserved for the ``gluing indices'' at the nodes.)
For simplicity of the notation, we adopt the \emph{summation convention} to
\emph{sum over all repeated indices}.


Let {$\displaystyle \frac{d \tau_G}{d \epsilon_s} := \hat{s}(z) \tau_G$}. Then
\[
 \begin{split}
  \frac{d G_0(\epsilon_s)}{d \epsilon_s} =
   &\sum_{l=1}^{\infty} \sum_{n=0}^{\infty}
        \sum_{i,j} (s_{l})_{ij} q^j_{l+n} \p_{q^i_n} G_0
   + \frac{1}{2} (-1)^n (s_{l+n+1})_{ij} q^i_n q^j_l . \\
  \frac{d G_g(\epsilon_s)}{d \epsilon_s} =
   &\sum_{l=1}^{\infty} \sum_{n=0}^{\infty}
        \sum_{i,j} (s_{l})_{ij} q^j_{l+n} \p_{q^i_n} G_g,
    \quad \text{for $g \ge 1$}. \\
 \end{split}
\]

Define
\[
  \la \p^{i_1}_{k_1} \p^{i_2}_{k_2} \ldots \p^{i_n}_{k_n} \ra_g
  := \frac{\p^n G_g}{\p \bt^{i_1}_{k_1} \p \bt^{i_2}_{k_2} \ldots
    \p \bt^{i_n}_{k_n}},
\]
and denote $\la \ldots \ra := \la \ldots \ra_0$. These functions $\la \ldots
\ra_g$ will be called \emph{axiomatic Gromov--Witten invariants}. Then
\begin{equation} \label{e:s0}
 \begin{split}
 &\frac{d}{d \epsilon_s} \langle \p^{i_1}_{k_1} \p^{i_2}_{k_2} \ldots \rangle\\
 = &\sum (s_l)_{ij} q^j_{l+n} \langle \p^i_n \p^{i_1}_{k_1} \ldots \rangle
  +\sum_{l=1}^{\infty} \sum_{i,a} (s_l)_{i i_a} \langle \p^i_{k_a-l}
   \p^{i_1}_{k_1} \ldots \hat{\p^{i_a}_{k_a}} \ldots \rangle \\
  + &\frac{\delta}{2} \Big( (-1)^{k_1} \sum (s_{k_1+k_2+1})_{i_1 i_2}
   + (-1)^{k_2} \sum (s_{k_1+k_2+1})_{i_2 i_1} \Big),
 \end{split}
\end{equation}
where $\delta=0$ when there are more than 2 insertions and
$\delta=1$ when there are two insertions. The notation
$\hat{\p^{i}_{k}}$ means that $\p^{i}_{k}$ is omitted from the summation.
We \emph{assume that there are at least two insertions},
as this is the case in our application.

For $g\ge 1$
\begin{equation} \label{e:sg}
 \begin{split}
  &\frac{d}{d\epsilon_s} \langle \p^{i_1}_{k_1} \p^{i_2}_{k_2} \ldots \rangle_g\\
  = &\sum (s_l)_{ij} q^j_{l+n} \langle \p^i_n \p^{i_1}_{k_1} \ldots \rangle_g
  + \sum \sum_{a} (s_l)_{i i_a} \langle \p^i_{k_a-l}
   \p^{i_1}_{k_1} \ldots \hat{\p^{i_a}_{k_a}} \ldots \rangle_g
 \end{split}
\end{equation}

\subsection{Quantization of upper triangular subgroups}

The quantization of $r(z)$ is
\[
 \begin{split}
 \hat{r}(z) = &\sum_{l=1}^{\infty} \sum_{n=0}^{\infty}
    \sum_{i,j} (r_l)_{ij} q^j_n \p_{q^i_{n+l}} \\
  + &\frac{\hbar}{2} \sum_{l=1}^{\infty} \sum_{m=0}^{l-1}
   (-1)^{m+1} \sum_{i j} (r_l)_{ij} \p_{q^i_{l-1-m}} \p_{q^j_m}.
 \end{split}
\]

%
%

Therefore
\begin{equation} \label{e:r0}
 \begin{split}
  \frac{d}{d \epsilon_r} &\langle \p^{i_1}_{k_1} \p^{i_2}_{k_2} \ldots \rangle
 = \sum_{l=1}^{\infty} \sum_{n=0}^{\infty} \sum_{i,j} (r_l)_{ij} q^j_n
   \langle \p^i_{n+l} \p^{i_1}_{k_1} \ldots \rangle \\
  + &\sum_{l=1}^{\infty} \sum_{i,a} (r_l)_{i i_a} \langle \p^i_{k_a+l}
   \p^{i_1}_{k_1} \ldots \hat{\p^{i_a}_{k_a}} \ldots \rangle \\
  + &\frac{1}{2} \sum_{l=1}^{\infty} \sum_{m=0}^{l-1} (-1)^{m+1}
   \sum_{i j} (r_l)_{ij} \p^{i_1}_{k_1} \p^{i_2}_{k_2} \ldots
   ( \langle \p^i_{l-1-m} \rangle \langle \p^j_m \rangle ) .
 \end{split}
\end{equation}

For $g \ge 1$
\begin{equation} \label{e:rg}
 \begin{split}
  &\frac{d \langle \p^{i_1}_{k_1} \p^{i_2}_{k_2} \ldots \rangle_g}
        {d \epsilon_r} \\
  = &\sum_{l=1}^{\infty} \sum_{n=0}^{\infty} \sum_{i,j} (r_l)_{ij} q^j_n
   \langle \p^i_{n+l} \p^{i_1}_{k_1} \ldots \rangle_g \\
  + &\sum_{l=1}^{\infty} \sum_{i,a} (r_l)_{i i_a} \langle \p^i_{k_a+l}
    \p^{i_1}_{k_1} \ldots \hat{\p^{i_a}_{k_a}} \ldots \rangle_g \\
  + &\frac{1}{2} \sum_{l=1}^{\infty} \sum_{m=0}^{l-1} (-1)^{m+1}
   \sum_{i j} (r_l)_{ij} \langle \p^i_{l-1-m} \p^j_m
    \p^{i_1}_{k_1} \p^{i_2}_{k_2} \ldots \rangle_{g-1} \\
  + &\frac{1}{2} \sum_{l=1}^{\infty} \sum_{m=0}^{l-1} (-1)^{m+1}
   \sum_{i j} \sum_{g'=0}^g (r_l)_{ij} \p^{i_1}_{k_1} \p^{i_2}_{k_2} \ldots
   ( \langle \p^i_{l-1-m} \rangle_{g'} \langle \p^j_m \rangle_{g-g'} ).
 \end{split}
\end{equation}

\section{Invariance of tautological equations}

\subsection{Invariance under lower triangular subgroups}

\begin{theorem} \label{t:sinv} \emph{($S$-invariance theorem)}
All tautological equations are invariant under action of lower triangular
subgroups of the twisted loop groups.
\end{theorem}

\begin{proof}
Let $E=0$ be a tautological equation of axiomatic Gromov--Witten invariants.
Suppose that this equation holds for a given semisimple Frobenius manifold,
e.g.~$H^{N pt} \cong \cc^N$. We will show that $\hat{s} E =0$. This will prove
the theorem.

$\hat{s} E =0$ follows from the following facts:
\begin{enumerate}
\item[(a)] The combined effect of the first term in \eqref{e:s0} (for genus
zero invariants) and in \eqref{e:sg} (for $g \ge 1$ invariants) vanishes.
\item[(b)] The combined effect of the remaining terms in \eqref{e:s0} and in
\eqref{e:sg} also vanishes.
\end{enumerate}
(a) is due to the fact that the sum of the contributions from the first
term is a derivative of the original equation $E=0$ with respect to $q$
variables. Therefore it vanishes.

It takes a little more work to show (b).
Recall that all tautological equations are induced from moduli spaces of
curves. Therefore, any relations of tautological classes on $\mbar_{g,n}$
contain no genus zero components of two or less marked points.
However, when one writes down the induced equation for
(axiomatic) Gromov--Witten invariants, the genus zero invariants with two
insertions will appear. This is due to the difference between the cotangent
classes on $\mbar_{g,n+m}(X, \beta)$ and the pull-back classes from
$\mbar_{g,n}$. Therefore the only contribution from the third term of
\eqref{e:s0} comes from these terms. More precisely, let $\psi_j$ (descendents)
denote the $j$-th cotangent class on $\mbar_{g,n+m}(X,\beta)$ and
$\bar{\psi}_j$ (ancestors) the pull-backs of cotangent classes from
$\mbar_{g,n}$ by the combination of the stabilization and forgetful morphisms
(forgetting the maps and extra marked points, and stabilizing if necessary).
Let $D_j$ be the divisor on
$\mbar_{g,n+m}(X,\beta)$ defined by the image of the gluing morphism
\[
 \sum_{\beta'+\beta''=\beta} \sum_{m' + m''=m}
 \mbar^{(j)}_{0, 2+m'}(X,\beta') \times_X \mbar_{g,n+m''}(X,\beta'')
 \to \mbar_{g,n+m}(X,\beta),
\]
where $\mbar_{g,n+m''}(X,\beta'')$ carries all first $n$ marked points
except the $j$-th one, which is carried by $\mbar^{(j)}_{0, 2+m'}(X,\beta')$.
It is easy to see geometrically that $\psi_j - \bar{\psi}_j = D_j$.
(See e.g.~\cite{KM}.) Let us denote $\la \p^{\mu}_{k,\bar{l}}, \ldots \ra$
the generalized (axiomatic) Gromov--Witten invariants with
$\psi_1^k \bar{\psi}_1^l \ev_1^*(\phi_{\mu})$ at the first marked point.
The above relation can be rephrased in terms of invariants as
\begin{equation*} 
  \la \p^i_{k,\bar{l}} \ldots \ra_g =
  \la \p^i_{k+1, \overline{l-1}} \ldots \ra_g
  - \la \p^i_k \p^{\mu} \ra \la \p^{\mu}_{\overline{l-1}} \ldots \ra_g.
\end{equation*}
Repeat this process of reducing $\bar{l}$, one can show by induction that
\[
 \begin{split}
  &\la \p^i_{k,\bar{l}} \ldots \ra_g =
  \la \p^i_{k+r, \overline{l-r}} \ldots \ra_g
  -\la \p^i_{k+r-1} \p^{\mu_1} \ra \la \p^{\mu_1}_{\overline{l-r}} \ldots\ra_g
  - \ldots \\
  &-\la \p^i_{k} \p^{\mu_1} \ra \left[ \sum_{p=1}^r (-1)^{p+1}
  \sum_{k_1 + \ldots + k_p = r-p} \la \p^{\mu_1}_{k_1} \p^{\mu_2} \ra
  \ldots \la \p^{\mu_{p-1}}_{k_{p-1}} \p^{\mu_p} \ra
  \la \p^{\mu_p}_{k_p, \overline{l-r}} \ldots \ra_g \right] .
 \end{split}
\]

Now suppose that one has an equation of tautological classes of $\mbar_{g,n}$.
Use the above equation (for $r=l$) one can translate the equation of
tautological classes on $\mbar_{g,n}$ into an equation of the (axiomatic)
Gromov--Witten invariants. The term-wise cancellation of the contributions
from the second and the third terms of \eqref{e:s0} and \eqref{e:sg} can
be seen easily by straightforward computation.
\end{proof}

If the above description is a bit abstract, the reader might want to
try the following simple example.
$\psi_1^2$ on $\mbar_{g,1}$ is translated into invariants:
\[
 \la \p^x_2 \ra_g - \la \p_1^x \p^{\mu} \ra \la \p^{\mu} \ra_g
  - \la \p^x \p^{\mu} \ra \la \p_1^{\mu} \ra_g
  + \la \p^x \p^{\mu} \ra \la \p^{\mu} \p^{\nu} \ra \la \p^{\nu} \ra_g .
\]
The above ``translation'' from tautological classes to Gromov--Witten
invariants are worked out explicitly in some examples in Sections~6 and 7 of
\cite{EG1}.

\subsection{Reduction to $q_0=0$}
The arguments in this section are mostly taken from \cite{GL}.

Let $E=0$ be a tautological equation of (axiomatic) Gromov--Witten invariants.
Since we have already proved $\hat{s}(E)=0$, our next goal would be
to show $\hat{r}(E)=0$. In this section, we will show that it suffices to
check $\hat{r}(E)=0$ on the subspace $\bq_0=0$.

\begin{lemma}
It suffices to show $\hat{r} E=0$ on each level set of the map $\bq \mapsto s$
in \eqref{e:dg}.
\end{lemma}

\begin{proof}
The union of the level sets is equal to $\HH_+$.
\end{proof}

\begin{lemma} \label{l:3}
It suffices to check the relation for all $\hat{r}(z) E =0$ along
$z\HH_{+}$ (i.e.~$\bq_0=0$).
\end{lemma}

\begin{proof}
It is proved in \cite{AG4} that a particular lower triangular matrix
$S_s$, which is called ``calibration'' of the Frobenius manifold, transforms
the level set at $s$ to $z \HH_+$. $S$-invariance Theorem then concludes the
proof.
\end{proof}

In fact, $S_s$ is a fundamental solution of the horizontal sections of the
Dubrovin (flat) connection, in $z^{-1}$ formal series. It was discovered in
\cite{AG2}, with preceding work in \cite{KM} and \cite{EG2}, that $\A:=
\hat{S}_s \tau^X$ is the corresponding generating function for ``ancestors''.
Therefore the transformed equation $\hat{S}_s E \hat{S}_s^{-1} =0$ is really an
equation of \emph{ancestors}.

\subsection{Invariance under upper triangular subgroup}

\begin{theorem} \emph{($R$-invariance theorem)}
The union of the sets of genus $g'$ equations for $g' \le g$ is invariant under
the action of upper triangular subgroup, for $g \le 2$. \footnote{We really
have not defined what we mean by ``the set of genus g equations'' for $g >2$,
as the current paper only deals with $g \le 2$ cases. One can think that it is
the set of all tautological equations at genus $g$.}
\end{theorem}

In fact, a stronger ``filtered'' statement holds. We will state the genus two
part:
\begin{enumerate}
\item[(I)] The combination of genus zero equations, genus one equations and
Mumford's equation is $R$-invariant.

\item[(II)] The combination of genus zero equations, genus one equations and
genus two Mumford's and Getzler's equations is $R$-invariant.

\item[(III)] The combination of genus zero equations, genus one equations and
genus two equations by Mumford, Getzler and BP is $R$-invariant.
\end{enumerate}

\begin{remark*}
1. $R$-invariance in $g=1$ is proved in \cite{GL}.

2. There are other genus two tautological equations, like the 6-point equation
discovered by Faber--Pandharipande (private communication). However, its role
in semisimple Gromov--Witten theory is not clear at this point.

3. $R$-invariance theorem is expected to hold for all $g$. In fact, under
plausible assumptions, $R$-invariance technique can be used to ``derive'' all
known tautological equations, including all tautological equations appeared
in this paper, and other new equations in higher genus.
This will be discussed in a separate paper.
\end{remark*}

The rest of the section is devoted to the proof of $R$-invariance theorem, and
therefore the Main Theorem. In fact, the proof of (I), (II) and (III) follow
the same line of arguments, so we will only treat (I) in details.

Recall that Mumford's genus two equation is of the form, \emph{in the
orthonormal basis} as usual, with summation convention,

\begin{equation} \label{e:m}
 \begin{split}
  M := &-\la \p^x_2 \ra_2 + \la \p_1^x \p^{\mu} \ra \la \p^{\mu} \ra_2
  + \la \p^x \p^{\mu} \ra \la \p_1^{\mu} \ra_2 \\
  &- \la \p^x \p^{\mu} \ra \la \p^{\mu} \p^{\nu} \ra \la \p^{\nu} \ra_2
  + \frac{7}{10} \la \p^x \p^{\mu} \p^{\nu} \ra
   \la \p^{\mu} \ra_1 \la \p^{\nu} \ra_1 \\
  &+ \frac{1}{10} \la \p^x \p^{\mu} \p^{\nu} \ra \la \p^{\mu} \p^{\nu} \ra_1
  -\frac{1}{240} \la \p^{\mu} \p^{\nu} \p^{\nu} \ra \la \p^x \p^{\mu} \ra_1 \\
  &+ \frac{13}{240} \la \p^x \p^{\mu} \p^{\mu} \p^{\nu} \ra \la \p^{\nu} \ra_1
  + \frac{1}{960} \la \p^x \p^{\mu} \p^{\mu} \p^{\nu} \p^{\nu} \ra =0 .
 \end{split}
\end{equation}

\begin{lemma} \label{l:l}
It suffices to check $(r(z)) \hat{}\ M =0$ for $l=1$ and $l=2$
(and on $q_0=0$).
\end{lemma}

\begin{proof}
It is easy to see that when $l \ge 3$, all terms in $(r(z)) \hat{} \ M$
\eqref{e:rm} vanish by the $(3g-2)$-jet property, which is satisfied for
the ancestor invariants (\cite{EG2} and Section~5 of \cite{AG2}).

In geometric terms, this is due to the fact that Mumford's equation \eqref{e:m}
is a codimension 2 tautological relation in $\mbar_{2,1}$, whose dimension
is equal to 4. Since $(r_l z^l) \hat{}$ carries codimension $k$ strata to
codimension $k+l$ ones, $(r_l z^l) \hat{} M =0$ for $l \ge 3$.
\end{proof}

In fact, as we will see, the checking of invariance is really straightforward
for $l=1$, and almost trivial for $l =2$. Note we will use
the following conventions:
\begin{itemize}
\item If some term does not contain $l$,
that means $l=1$ and it vanishes for $l \ge 2$.
\item $\p^x$ is any flat vector field, e.g.~$\p^x=\p^{\mu}_k$. Therefore
$\p^x_2$ means the descendent index is at least two, but could be greater.
\end{itemize}

\begin{equation} \label{e:rm}
 \begin{split}
  (r(z)) \hat{} \ M & \\
   = \sum_l (r_l)_{ij} \bigg[
   &-\frac{1}{2} \sum (-1)^{m+1} \la \p_{l-1+m}^i \p^j_m \p^x_2 \ra_1 \\
  &- \sum (-1)^{m+1} \la \p^i_{l-1-m} \p^x_2 \ra_1 \la \p^j_m \ra_1 \\
  &+\frac{7}{5} \sum \la \p^x \p^j \p^{\nu} \ra \la \p^i_1 \ra_1
    \la \p^{\nu} \ra_1 \\
  & -\frac{7}{10} \sum \la \p^x \p^{\mu} \p^{\nu} \ra
    \la \p^i \p^j \p^{\mu} \ra \la \p^{\nu} \ra_1 \\
  &+\frac{1}{5} \sum \la \p^x \p^j \p^{\nu} \ra \la \p^i_l \p^{\nu} \ra_1 \\
  &+\frac{1}{20} \sum (-1)^{m+1} \la \p^x \p^{\mu} \p^{\nu} \ra
   \la \p^i_{l-1-m} \p^j_m \p^{\mu} \p^{\nu} \ra \\
  &+\frac{1}{10} \sum (-1)^l \la \p^x \p^{\mu} \p^{\nu} \ra
   \la \p^i \p^{\mu} \p^{\nu} \ra \la \p^j_{l-1} \ra_1 \\
  &-\frac{1}{240} \sum \la \p^j \p^{\nu} \p^{\nu} \ra \la \p^x \p^i_l \ra_1 \\
  &- \frac{1}{480} \sum (-1)^{m+1} \la \p^{\mu} \p^{\nu} \p^{\nu} \ra
  \la \p^i_{l-1-m} \p^j_m \p^x \p^{\mu} \ra \\
  &- \frac{1}{240} \sum (-1)^l \la \p^{\mu} \p^{\nu} \p^{\nu} \ra
    \la \p^i \p^x \p^{\mu} \ra \la \p^j_{l-1} \ra_1  \\
  &+\frac{13}{120} \sum \la \p^x \p^i_l \p^j \p^{\nu} \ra \la \p^{\nu} \ra_1\\
  &+\frac{13}{240} \sum \la \p^x \p^{\mu} \p^{\mu} \p^i_l \ra \la \p^j \ra_1 \\
  &+\frac{13}{240} \sum \la \p^x \p^{\mu} \p^{\mu} \p^j \ra \la \p^i_l \ra_1 \\
  &-\frac{13}{480} \sum \la \p^x \p^{\mu} \p^{\mu} \p^{\nu} \ra
  \la \p^i \p^j \p^{\nu} \ra \\
  &-\frac{13}{120} \sum \la \p^i \p^x \p^{\mu} \ra
  \la \p^j \p^{\mu} \p^{\nu} \ra \la \p^{\nu} \ra_1 \\
  &-\frac{13}{240} \sum \la \p^i \p^x \p^{\nu} \ra
  \la \p^j \p^{\mu} \p^{\mu} \ra \la \p^{\nu} \ra_1 \\
  &+ \frac{1}{240} \sum \la \p^x \p^i_l \p^j \p^{\nu} \p^{\nu} \ra \\
  &- \frac{1}{480} \sum \la \p^i_{l-1} \p^x \p^{\mu} \p^{\mu} \ra
  \la \p^j \p^{\nu} \p^{\nu} \ra \\
  &- \frac{1}{240} \sum \la \p^i_{l-1} \p^x \p^{\mu} \p^{\nu} \ra
  \la \p^j \p^{\mu} \p^{\nu} \ra \\
  &- \frac{1}{240} \sum \la \p^i_{l-1} \p^x \p^{\mu} \ra
  \la \p^j \p^{\mu} \p^{\nu} \p^{\nu} \ra \bigg]
 \end{split}
\end{equation}

In the above, the following observations simplify the calculation:
\begin{itemize}
\item The contributions from first terms of \eqref{e:r0} \eqref{e:rg} vanish.
\item The contributions from the second terms of \eqref{e:r0} \eqref{e:rg},
when acting on $\p^x_k$, cancell with each other.
\item The contributions from the last terms of \eqref{e:r0} \eqref{e:rg},
when all $\p^{i_a}_{k_a}$ act on one of the factor with highest genus
in a term, cancell each other.
\end{itemize}
The first is already explained in the proof of Theorem~\ref{t:sinv}. The second
is due to the fact that the sum of the second term merely changes one flat
vector to another: $\p^{\mu}_k \mapsto \sum (r_l)_{i\mu} \p^i_{k+l}$.
Since $\p^x$ is an arbitrary flat vector field, this change sums to zero.
The last case is due to the fact that it sums to a derivative of $M=0$ with
respect to $\bq$.

The two cases of $l=1$ and $l=2$ will be discussed separately.

\subsubsection{The case $l=1$}
Since there are only invariants in genus zero and one, all descendents can
be removed by TRR. (The flat vector field $\p^x$ might contain
descendents, but we will only remove the ``apparent descendents''.)

After applying TRRs, all terms can be classified into two types:
\begin{itemize}
\item Type 1: terms involving genus one invariants and
\item Type 2: terms with only genus zero invariants.
\end{itemize}

\textbf{Type 1:} After obvious vanishing (on the dimensional ground), the
contribution is
\[
 \begin{split}
  \sum (r_1)_{ij} \bigg[
  &\frac{1}{2} \la \p^x_1 \p^i \p^j \p^{\mu} \ra \la \p^{\mu} \ra_1
  + \frac{1}{24} \la \p^x_1 \p^i \p^{\mu} \p^{\mu} \ra \la \p^j \ra_1 \\
  +&\frac{7}{120} \la \p^i \p^{\mu} \p^{\mu} \ra \la \p^x \p^j \p^{\nu} \ra
  \la \p^{\nu} \ra_1
  -\frac{7}{10} \la \p^x \p^{\mu} \p^{\nu} \ra \la \p^i \p^j \p^{\mu} \ra
  \la \p^{\nu} \ra_1  \\
  +&\frac{1}{5}  \la \p^x \p^j \p^{\nu} \ra \la \p^i \p^{\mu} \p^{\nu} \ra
  \la \p^{\mu} \ra_1
  -\frac{1}{10}  \la \p^x \p^{\mu} \p^{\nu} \ra \la \p^i \p^{\mu} \p^{\nu} \ra
  \la \p^j \ra_1 \\
  -&\frac{1}{240} \la \p^x \p^i \p^{\mu} \ra \la \p^j \p^{\nu} \p^{\nu} \ra
  \la \p^{\mu} \ra_1
  +\frac{1}{240} \la \p^x \p^i \p^{\mu} \ra \la \p^{\mu} \p^{\nu} \p^{\nu} \ra
  \la \p^j \ra_1 \\
  +&\frac{13}{240} \la \p^x \p^i_1 \p^j \p^{\mu} \ra \la \p^{\mu} \ra_1
  +\frac{13}{240} \la \p^x \p^i_1 \p^{\mu} \p^{\mu} \ra \la \p^j \ra_1 \\
  -&\frac{13}{120} \la \p^x \p^i \p^{\mu} \ra \la \p^j \p^{\mu} \p^{\nu} \ra
  \la \p^{\nu} \ra_1
  -\frac{13}{120} \la \p^x \p^i \p^{\mu} \ra \la \p^j \p^{\nu} \p^{\nu} \ra
  \la \p^{\nu} \ra_1 \bigg]
 \end{split}
\]

The above sum can be put into two groups. One with genus one invariants
of the form $\la \p^j \ra_1$, the other with $\la \p^{\mu} \ra_1$.
It is easy to see that each group gets exact cancellation by genus zero
TRR (and WDVV equation, which is a consequence of TRR).

\textbf{Type 2:}
\[
 \begin{split}
  \sum (r_1)_{ij} \bigg[
  &\frac{1}{48} \la \p^x_1 \p^i \p^j \p^{\mu} \p^{\mu} \ra
  +\frac{1}{120} \la \p^x \p^j \p^{\nu} \ra
  \la \p^i \p^{\nu} \p^{\mu} \p^{\mu} \ra \\
  -&\frac{1}{20} \la \p^i \p^j  \p^{\mu} \p^{\nu} \ra
  \la \p^x \p^{\nu} \p^{\mu} \ra
  -\frac{1}{5760} \la \p^i \p^x \p^{\mu} \p^{\mu} \ra
  \la \p^j \p^{\nu} \p^{\nu} \ra \\
  +&\frac{1}{480} \la \p^i \p^j \p^x \p^{\mu} \ra
  \la \p^{\mu} \p^{\nu} \p^{\nu} \ra
  +\frac{13}{5760} \la \p^x \p^j \p^{\mu} \p^{\mu} \ra
  \la \p^i \p^{\nu} \p^{\nu} \ra \\
  -&\frac{13}{480} \la \p^x \p^{\mu} \p^{\mu} \p^{\nu} \ra
  \la \p^i \p^j \p^{\nu} \ra
  +\frac{1}{240} \la \p^x \p^{\mu} \p^{\mu} \p^i_1 \p^j \ra \\
  +&\frac{1}{240} \la \p^x \p^i \p^{\mu} \ra
  \la \p^j \p^{\mu} \p^{\nu} \p^{\nu} \ra
  -\frac{1}{480} \la \p^x \p^i \p^{\mu} \p^{\mu} \ra
  \la \p^j \p^{\nu} \p^{\nu} \ra \\
  -&\frac{1}{240} \la \p^x \p^i \p^{\mu} \p^{\nu} \ra
  \la \p^j \p^{\mu} \p^{\nu} \ra \bigg]
  \end{split}
\]

This can be seen to vanish by elementary calculation, utilizing
genus zero TRR and WDVV.

\subsubsection{The case $l=2$}
Due to \eqref{e:r}, $(r_2)_{ij} = - (r_2)_{ji}$. However, all terms in
\eqref{e:rm} with nonvanishing contribution at $l=2$ are all \emph{symmetric}
in $i$ and $j$. For example the first term contributes
\[
 \begin{split}
  & -\frac{1}{2} \sum_{i j} (r_2)_{ij} \la \p_1^i \p^j \p^x_2 \ra_1
   +\frac{1}{2} \sum_{i j} (r_2)_{ij} \la \p^i \p^j_1 \p^x_2 \ra_1 \\
  = &\sum_{i j} (r_2)_{ij} \la \p^i \p^j_1 \p^x_2 \ra_1
  = \sum_{i j} (r_2)_{ij} [\la \p^i \p^j \p^{\mu} \ra \la \p^x_2 \p^{\mu} \ra
  + \frac{1}{24} \la \p^i \p^j \p^{\mu} \p^{\mu} \p^x_2 \ra ],
 \end{split}
\]
which vanishes as $(r_2)_{ij}$ in anti-symmetric while the other factor is
symmetric in $i$ and $j$.

The terms do not cancell by itself after summation over $i,j$ will contribute
zero by cancellation. More precisely, the summation is
\[
 \sum (r_2)_{ij}
 (\frac{-1}{576} +\frac{1}{240} -\frac{1}{5760} -\frac{1}{5760}-\frac{1}{480})
 \la \p^i \p^x \p^{\mu} \ra \la \p^{\mu} \p^{nu} \p^{\nu} \ra
 \la \p^j \p^{\alpha} \p^{\alpha} \ra ,
\]
which is easily seen to vanish.

\subsubsection{$R$-invariance of Getzler's and BP's equations}

The $R$-invariance of Getzler's and BP's equations are proved in the same way.
It is not very illuminating to write down all the calculations here as the
proof follows the same line of arguments and the formula is quite lengthy.
Therefore, we will list only some main steps here.
Recall that $q_0=0$ is always assumed.

For Getzler's genus two equation:
\begin{enumerate}
\item $G':= (r(z)) \hat{} \ (\text{Getzler's equation})$ involves only up to
$l=3$, using the same argument in Lemma~\ref{l:l}.

\item $G'$ contains only one term with genus two component:
$3 \la \p^x \p^y \p^{\mu} \ra \la \p^{\mu}_{l+1} \ra_2$. This can be written,
via Mumford's equation as a summation of genus 1 and 0 invariants only.

\item One may remove the ``apparent descendents'' in genus one invariants by
TRR. Here the apparent descendents means the lower indices in our notation.
Although $\p^x$ might already contain descendents implicitly, they won't be
removed.

\item One can group terms together according to the types of the factors of the
genus one invariants in these terms. Since there is no additional relation in
genus one, they will cancell within each group.
We note that no 4-point genus one invariants appeared in the calculation.
Therefore Getzler's genus one equation is actually never used.

\item The cancellations of the genus zero part of \eqref{e:rm} involves
only (derivatives of) WDVV equation.
\end{enumerate}
There are only two different points in the proof of $R$-invariance of BP's
equation. First, the calculation will go up to $l=4$. Second,
Getzler's genus two equation will be used.


\end{document}